\documentclass[11pt]{amsart}

\issueinfo{00}{0}{Month}{Year}
\copyrightinfo{2008}{University of Calgary}
\pagespan{1}{2}
\PII{ISSN 1715-0868}
\include{cdmstdcmds}

\usepackage{amsrefs}	
\usepackage{amsfonts}
\usepackage{amsmath, amsthm, amssymb}
\usepackage[capitalize, noabbrev]{cleveref}
\newtheorem{theorem}{Theorem}[section]
\newtheorem{corollary}[theorem]{Corollary}
\newenvironment{pf}           {\noindent{\bf Proof:} }%
                                {\null\hfill$\Box$\par\medskip}

\begin{document}

\title{Total colourings of direct product graphs}

\author{Kyle MacKeigan}
\author{Jeannette Janssen}
\address{}
\email{kyle.m.mackeigan@gmail.com}
\email{jeannette.janssen@dal.ca}
%

\date{(date1), and in revised form (date2).}
\subjclass[2000]{05C15, 05C76}
\keywords{Graph theory, graph products, total colouring, combinatorics}

\thanks{This research was supported by an NSERC Discovery grant}

\begin{abstract}
A graph is $k$-total colourable if there is an assignment of $k$ different colours to the vertices and edges of the graph such that no two adjacent nor incident elements receive the same colour. The total chromatic number of some direct product graphs are determined. In particular, a sufficient condition is given for direct products of bipartite graphs to have total chromatic number equal to its maximum degree plus one. Partial results towards the total chromatic number of $K_n\times K_m$ are also established.
\end{abstract}

\maketitle

\section{Introduction}

All graphs considered in this paper are finite and simple. A total colouring is an assignment of colours to vertices and edges (elements), so that any two adjacent or incident elements receive different colours. The minimum number of colours required for a total colouring is called the total chromatic number of $G$, denoted $\chi''(G)$. Bezhad \cites{Bezhad1,Bezhad2} and Vizing \cite{Vizing} independently conjectured, in what is now known as the Total Colouring Conjecture, that for any graph $G$, the following inequality holds:
$$\Delta(G)+1\leq \chi''(G)\leq\Delta(G)+2$$
where $\Delta(G)$ is the maximum degree of $G$. Graphs that require only $\Delta(G)+1$ colours are called type I graphs whereas graphs that require  $\Delta(G)+2$ colours are called type II graphs.

The direct product of two graphs $G$ and $H$ is a graph, denoted $G\times H$, whose vertex set is $V(G)\times V(H)$ and for which vertices $(u,v)$ and $(u',v')$ are adjacent if and only if $uu'\in E(G)$ and $vv' \in E(H)$. The following work was done on determining the total chromatic number of $K_n\times K_m$. Yap \cites{Yap1,Yap2} showed that for $n>2$, $K_n\times K_2$ is a type I graph. For $n=2$,   $K_2\times K_2\simeq 2K_2$ is a type II graph. Geetha and Somasundaram \cite{Geetha} proved that if $n\geq 4$ and $n$ is even, then $K_n\times K_n$ is a type I graph.

In this paper, it is shown that if $G\times K_2$ is type I, then $G\times H$ is type I for any bipartite graph $H$. It is then shown that if $n$ or $m$ is even and $n,m\geq 3$, then $K_n\times K_m$ is a type I graph. 
 
\section{Direct Product With Complete Graphs.}

For the Cartesian graph product, Kemnitz and Marangio \cite{Kemnitz} proved that if $G\square K_2$ is type I, then $G\square H$ is type I where $H$ is any bipartite graph. The analogous result has also been proved for the strong graph product \cite{Geetha}. Here the result for the direct graph product is established. The Cartesian and strong graph product proofs utilized that the product graph had copies of $G$ and $H$. This is not always the case for the direct product. However, the direct product with a bipartite graph is a bipartite graph, and it is well known that bipartite graphs have $\Delta$-edge colourings.

\begin{theorem}\label{thrm: tensorproduct}
If $G\times K_2$ is type I then $G\times H$ is type I for any bipartite graph $H$.
\end{theorem}
\begin{pf} If $H$ is edgeless, then $G\times H$ is edgeless and thus a type I graph. Suppose now that $H$ is a bipartite graph with $\Delta(H)\geq 1$ and the following bipartition of the vertices: $\{x_1,x_2,\dots, x_n\}$ and $\{y_1,y_2,\dots,y_m\}$ where $n\leq m$. Suppose that $G\times K_2$ is a type I graph where $G$ has vertices $v_1,v_2,\dots, v_r$ and $K_2$ has vertices $z_1$ and $z_2$. For $1\leq i \leq n$ and $1\leq j\leq m$, partition the vertices of $G\times H$ into the following sets: $X_i=\{(v_k,x_i):1\leq k\leq r\}$ and $Y_j=\{(v_k,y_j):1\leq k \leq r\}$.

Let $f$ be a total colouring of $G\times K_2$. We now form a total colouring $g$ of $G\times H$. For each $v_k$ and each $x_i$, set $g((v_k,x_i))=f((v_k,z_1))$. Similarly, for each $v_k$ and each $y_j$, set $g((v_k,y_j))=f((v_k,z_2))$. Since $\{x_1,x_2,\dots,x_m\}$ and $\{y_1,y_2,\dots,y_m\}$ is a bipartition of $H$, none of the vertices in any of the sets $X_i$ are mutually adjacent. Similarly none of the vertices in $Y_j$ are mutually adjacent. Now, suppose there are two adjacent vertices, $(v_k,x_i)$ and $(v_t,y_j)$. Then, by definition, $v_k$ and $v_t$ are adjacent in $G$, and thus, $(v_k,z_1)$ and $(v_t,z_2)$ are adjacent vertices in $G\times K_2$, and thus receive different colours in $f$. Therefore, there is no conflict between vertices. It remains to extend this colouring to the edges.

Now, since $H$ is a bipartite graph, it is $\Delta(H)$-edge colourable. Take any one edge colour class $C$ in an edge colouring of $H$. Since this is a $\Delta(H)$-edge colouring, all of the vertices with maximum degree in $H$ will be an endpoint of one edge in $C$. If $x_iy_j\in C$, then in $G\times H$, set $g((v_k,x_i)(v_t,y_j))=f((v_k,z_1)(v_t,z_2))$. This is a proper edge colouring. Namely, suppose there are two incident edges, $(v_k,x_i)(v_t,y_j)$ and $(v_k,x_i)(v_{t'},y_{j'})$ with $x_iy_j\in C$ and $x_iy_{j'}\in C$. Then, by definition, $v_kv_t$ and $v_kv_{t'}$ are incident edges in $G$, and thus, $(v_k,z_1)(v_t,z_2)$ and $(v_k,z_1)(v_{t'},z_2)$ are incident edges in $G\times K_2$, and thus receive different colours in $f$. A similar argument holds for incident edges that share an endpoint $(v_t,y_j)$. Moreover, for each edge $(v_k,x_i)(v_t,y_j)$ with $x_iy_j\in C$, $g((v_k,x_i)(v_t,y_j))=f((v_k,z_1)(v_t,z_2)) \neq f((v_k,z_1))= g((v_k,x_i))$. Same argument shows that $g((v_k,x_i)(v_t,y_j))\neq g((v_t,y_j))$. Therefore, there is no conflict between vertices and edges.

Notice in $G\times H$ that $\deg((v,w))=\deg(v)\deg(w)$. Therefore, all vertices in $G\times H$ have degree that is at most a multiple of $\Delta(H)$. In particular, all vertices that are not of max degree, $\Delta(G)\Delta(H)$, have degree at most $\Delta(G)(\Delta(H)-1)$. Consider the graph obtained by removing all edges $(v_k,x_i)(v_t,y_j)$ with $x_iy_j\in C$ from $G\times H$. Since $C$ is incident with all vertices of maximum degree in $H$, the vertices in this graph will have maximum degree $\Delta(G)(\Delta(H)-1)$. But $G\times H$ is a bipartite graph with bipartition $X_1\cup X_2\cup \dots \cup X_n$ and $Y_1\cup Y_2\cup \dots \cup Y_m$. Thus the remaining edges can be coloured with $\Delta(G)(\Delta(H)-1)$ colours. Therefore $\Delta(G)(\Delta(H)-1)+ \Delta(G)+1 =\Delta(G)\Delta(H)+1=\Delta(G\times H)+1$ colours are used in this total colouring of $G\times H$. 
\end{pf}

Now one can determine that $G\times H$ is a type I graph, for any bipartite graph $H$, by showing that $G\times K_2$ is a type I graph. For instance, applying Theorem \ref{thrm: tensorproduct} to the fact that $K_n\times K_2$ is a type I graph for $n>2$ \cites{Yap1,Yap2}, the following corollary is established.

\begin{corollary}
If $n\neq 2$, then $K_n\times H$ is a type I graph for any bipartite graph $H$.
\end{corollary}

The graph $K_n\times K_2$ can be described as $K_{n,n}$ with a perfect matching removed. This graph is also sometimes described at the crown graph, denoted $J_{2n}$. Using a particular total colouring of $J_{2n}$, the total chromatic number of $K_n\times K_m$ when $n$ or $m$ is even and $n,m\geq 3$ is determined. 

\begin{theorem}\label{thrm:All The Crown}
If $n$ or $m$ is even and $n,m\geq 3$, then $K_n\times K_m$ is type I.
\end{theorem}
\begin{pf}
Suppose without loss of generality that $n\geq 4$ is even and label the vertices of $K_n$ as $v_1,v_2,\dots,v_n$ and the vertices of $K_m$ as $u_1,u_2,\dots, u_m$. Partition the vertices of $K_n\times K_m$ into the following sets:
\begin{align*}
X_1&=\{(v_1,u_k):1\leq k\leq m\}\\
X_2&= \{(v_2,u_k):1\leq k\leq m\}\\
&\vdots\\
X_{n}&= \{(v_n,u_k):1\leq k \leq m\}
\end{align*}

Note that the subgraph induced on $X_1\cup X_2$ is isomorphic to $J_{2m}$. The same is true for the subgraphs induced on $X_3\cup X_4,\dots,X_{n-1}\cup X_n$. The next step is to find a total colouring of $J_{2m}$ where the colours assigned to the vertices in each bipartition are distinct. Such a colouring can be acquired by using a particular edge colouring of $K_{m,m}$. 

A perfect rainbow matching in a graph is defined as a perfect matching, such that all the edges in the perfect matching have distinct colors. Yap \cite{Yap2} showed that for $m\geq 3$, there is an $m$ edge colouring of $K_{m,m}$ with a perfect rainbow matching. By removing the edges in this perfect rainbow matching from $K_{m,m}$ and colouring the endpoints of each edge by the colour of the removed edge, a total colouring of $J_{2m}$ is acquired. Let $h$ be this constructed total colouring. Note that the vertices in each bipartition are distinctively coloured. This is because each edge in the perfect rainbow matching were distinctively coloured.

Now, let $\{x_1,x_2,\dots ,x_m\}$ and $\{y_1,y_2,\dots ,y_m\}$ be a bipartition of $J_{2m}$, where $x_ky_k\not\in E(J_{2m})$. Then by the construction of $h$, $h(x_k)=h(x_t)$ and $h(y_k)=h(y_t)$ if and only if $k=t$. Also, $h(x_k)=h(y_t)$ if and only if $k=t$. We now form a total colouring $g$ of $K_n\times K_m$. For each $k$, $1\leq k\leq m$, and for each $i$, $1\leq i \leq n$, set $g((v_i,u_k))= h(x_k)$. Two vertices $(v_i,u_k)$ and $(v_j,u_t)$ are adjacent only if $k\neq t$, therefore $g((v_i,u_k))=h(x_k)\neq h(x_t)=g((v_j,u_t))$. Therefore, there is no conflict between the vertices. It remains to extend this colouring to the edges.

Consider an edge colouring $l$ of $K_n$ using colours $0,1,\dots,n-2$. Note that such a colouring exists because $n$ is even. If $l(v_iv_j)=0$, then for each $k$ and each $t$, $k\neq t$, $1\leq k,t \leq m$, set $g((v_i,u_k)(v_j,u_t))=h(x_k y_t)$. This is a proper edge colouring. Namely, suppose there are two incident edges, $(v_i,u_k)(v_j,u_t)$ and $(v_i,u_k)(v_{j'},u_{t'})$ where $l(v_iv_j)=l(v_iv_{j'})=0$. Then $t\neq t'$, $k\neq t$, and $k\neq t'$, therefore $x_ky_t$ and $x_ky_{t'}$ are incident edges in $J_{2m}$, and thus receive different colours in $h$. A similar argument holds for edges that share an endpoint $(v_j,u_t)$. Moreover, for each edge $(v_i,u_k)(v_j,u_t)$, where $l(v_iv_j)=0$, $g((v_i,u_k)(v_j,u_t))=h(x_ky_t)\neq h(x_k)= g((v_i,u_k))$ and $g((v_i,u_k)(v_j,u_t))=h(x_ky_t)\neq h(y_t)= g((v_j,u_t))$. Thus, there is no conflict between vertices and edges.

Now, consider an edge colouring $f$ of $J_{2m}$ with colours $0,1,\dots,m-2$. For $1\leq c\leq n-2$, if $l(v_iv_j)=c$, then for each $k$ and each $t$, $1\leq k,t\leq m$, $k\neq t$, set $g((v_i,u_k)(v_j,u_t))=c(m-1)+f(x_k y_t)+1$. These edges are assigned colours greater than $m-1$, so they will not conflict with the colours assigned to the vertices or the previously coloured edges. To show that these edges also do not conflict with each other, suppose that $(v_i,u_k)(v_j,u_t)$ and $(v_z,u_r)(v_j,u_t)$ are incident. If $i\neq z$, then $l(v_iv_j)=c_1$ and $l(v_zv_j)=c_2$ where $c_1\neq c_2$. Without loss of generality, suppose that $c_1<c_2$, then 
\begin{align*}
g((v_i,u_k)(v_j,u_t))&=c_1(m-1)+f(x_k y_t)+1\\
&\leq c_1(m-1)+m-2+1\\
&=c_1(m-1)+(m-1)\\
&=(m-1)(c_1+1)\\
&\leq (m-1)c_2\\
&<(m-1)c_2+f(x_ry_t)+1\\
&=g((v_z,u_r)(v_j,u_t))
\end{align*}

If $i=z$, then $g((v_i,u_k)(v_j,u_t))-g((v_i,u_r)(v_j,u_t))=f(x_ky_t)-f(x_ry_t)\neq 0$ because $x_ky_t$ and $x_ry_t$ are incident in $J_{2m}$ and $f$ is a proper edge colouring. A similar argument holds for incident edges that share an endpoint $(v_i,u_k)$. Therefore, $(n-2)(m-1)+ m=(n-1)(m-1)+1= \Delta(K_n\times K_m)+1$ colours are used in this total colouring of $K_n\times K_m$. 
\end{pf}

It is still an open problem to determine the total chromatic number of $K_n\times K_m$ when $n$ and $m$ are both odd.

\bibliographystyle{amsplain}
\bibliography{references}

\begin{bibdiv}
\begin{biblist}

\bib{Kemnitz}{article}{
      author={A.~Kemnitz, M.~Marangio},
       title={Total colorings of cartesian products of graphs},
        date={2003},
     journal={Congr. Numer.},
      volume={165},
       pages={99\ndash 109},
}

\bib{Bezhad1}{thesis}{
      author={Behzad, M.},
       title={Graphs and their chromatic number},
        type={Ph.D. Thesis},
     address={East Lansing},
        date={1965},
}

\bib{Bezhad2}{article}{
      author={Behzad, M.},
       title={The total chromatic number},
        date={1971},
     journal={Combinatorial Mathematics and Its Applications},
       pages={1\ndash 8},
}

\bib{Geetha}{article}{
      author={J.~Geetha, K.~Somasundaram},
       title={Total colorings of product graphs},
        date={2018},
     journal={Graphs and Combinatorics},
      volume={34},
       pages={339\ndash 347},
}

\bib{Vizing}{article}{
      author={Vizing, V.G.},
       title={Some unsolved problems in graph theory},
        date={1968},
     journal={Russian Mathematical Surveys},
      volume={23},
       pages={125\ndash 141},
}

\bib{Yap1}{article}{
      author={Yap, H.P.},
       title={Generalizations of two results of hilton on total colourings of a
  graph},
        date={1995},
     journal={Discrete Math.},
      volume={140},
       pages={245\ndash 252},
}

\bib{Yap2}{article}{
      author={Yap, H.P.},
       title={Total colourings of graphs},
        date={1996},
     journal={Lecture Notes in Mathematics 1623},
}

\end{biblist}
\end{bibdiv}

\end{document}